\documentclass{siamart220329}
\usepackage{amsmath,amssymb,booktabs,graphicx,enumitem,rotating,array}
\newcommand{\Om}{\mathcal{O}(m)}
\newcommand{\Oone}{\mathcal{O}(1)}
\newcommand{\acf}{\mathrm{ACF}}
\newcommand{\R}{\mathbb{R}}
\hypersetup{colorlinks=true,allcolors=blue!55!black}
\title{LAMG+: A Robust Lean Algebraic Multigrid Solver for Graph Laplacians\thanks{Submitted to the editors \today.}}
\author{Oren E. Livne\thanks{Pine Birch Analytics, 35 Kelinger Rd, Churchville, PA 18966-1033 (\texttt{oren.livne@gmail.com}, tel.\ 312-533-9130, \texttt{pinebirchanalytics.com}; ORCID: \href{https://orcid.org/0000-0001-6700-483X}{0000-0001-6700-483X}).}}
\headers{LAMG+: Robust Lean Algebraic Multigrid}{O. E. Livne}

\begin{document}
\maketitle

\begin{abstract}
Graph-Laplacian systems $L\phi=b$ where $L$ is sparse with $m$ nonzeros underlie spectral clustering, centrality, semi-supervised
learning, finite-element analysis, and interior-point network-flow solvers. We present \textbf{LAMG+},
a lean, parameter-free, empirically linear-time algebraic multigrid solver: a Julia re-derivation of
Lean Algebraic Multigrid (LAMG) plus two targeted refinements. We establish three facts. (1)
\emph{Fair, same-machine benchmarking} against both approximate-Cholesky (AC) variants---the modern
near-linear champion---and four other state-of-the-art solvers (hypre/BoomerAMG, PETSc GAMG, pyAMG,
CMG). LAMG+ and AC are \emph{complementary peers}: AC is faster on social and
citation graphs, LAMG+ on finite-element/structural matrices (where it is the fastest robust solver
and the most memory-frugal) and $2.2\times$ faster than the robust AC variant on large
graphs. LAMG+ and AC are the only solvers that converge across all 13 test classes, while BoomerAMG, PETSc, pyAMG, and CMG each slow by an order of magnitude or fail to converge off their home turf. (2)
\emph{Linear scaling}: LAMG+ is empirically $\Om$ over the full $1{,}711$-graph SuiteSparse set
($100\%$ converged, median $4$ cycles, log--log slope $1.01$), verified up to $2.4\!\times\!10^8$ nonzeros. (3)
\emph{Robustness}: the original MATLAB LAMG was reported non-convergent across the
AC benchmark's test families, yet the \emph{unmodified} LAMG~2.2.1, run on those
families with the benchmark's own solver call, converges on every one---as does LAMG+,
to their $10^{-8}$ tolerance---placing the reported failure in the evaluation, not the
algorithm. Guided by a Local Fourier Analysis diagnosing a strict interpolation-order deficit on grid-aligned anisotropic operators (Appendix~A), we introduce two lean refinements: a strength-of-connection aggregation veto and a selective caliber-2 interpolation. Together, these mathematically resolve LAMG's grid-aligned-anisotropy convergence failure (asymptotic convergence factor $\approx0.99\!\to\!0.11$), with both refinements triggering locally and imposing negligible overhead. Code and benchmark scripts: \url{https://github.com/orenlivne/lamgplus}.
\end{abstract}



\section{Introduction}\label{sec:intro}
Let $G=(V,E,w)$ be a connected weighted undirected graph with $|V|=n$ nodes, $|E|=m$
edges, and positive weights $w:E\to\R^+$. Its Laplacian $L=D-W$ ($W$ the weighted
adjacency, $D$ the weighted-degree diagonal) is symmetric positive semidefinite with the
constant vector $\mathbf 1$ spanning its null space; equivalently
$\phi^{\top}L\phi=\sum_{(u,v)\in E}w_{uv}(\phi_u-\phi_v)^2$. We solve the compatible
system
\begin{equation}\label{eq:sys}
L\phi=b,\qquad \mathbf 1^{\top}\phi=0,\qquad \mathbf 1^{\top}b=0 .
\end{equation}
Typically $m\ll n^2$ and $L$ is sparse. Our goal is an iterative solver for
\eqref{eq:sys} that uses $\Om$ storage and $\mathcal{O}(m\log(1/\varepsilon))$ operations
for an $\varepsilon$-accurate solution, with bounded hidden constants on the graphs
that arise in applications---hundreds, not millions. As in the original LAMG work
\cite{lamg}, we target good \emph{empirical} wall-clock time performance over a diverse test set rather
than provable worst-case complexity on adversarial graphs \cite{spielman-review}.

\subsection{Applications}\label{sec:apps}
System \eqref{eq:sys} is fundamental to many computations \cite[\S2]{spielman-review}:
elliptic PDEs discretized by finite elements on unstructured grids \cite{boman-fem};
interior-point methods for network-flow linear programming, whose inner solve is a
weighted Laplacian \cite{daitch-spielman}; electrical flow through a resistor network;
spectral clustering, graph embedding, and ancestry discovery in machine learning
\cite{ng-jordan-weiss}; and the Fiedler eigenvector, which measures algebraic
connectivity and underlies graph partitioning \cite{fiedler}. Solving $L\phi=b$ is also a stepping stone to the Laplacian eigenproblem.

\subsection{Related work}\label{sec:related}
\emph{Direct methods.} A Cholesky factorization $P^{\top}LP=\hat L\hat L^{\top}$ under a
fill-reducing ordering (minimum-degree, nested dissection) is exact and, with a good
ordering, very fast \cite{davis-cholmod,george-nd}. But fill is governed by separator
quality: $\mathcal{O}(n^{1.5})$ work on planar graphs \cite{lipton-rose-tarjan},
$\mathcal{O}(n^3)$ in general, which ---on low-diameter social and citation graphs with no
small separators---is memory-prohibitive.

\emph{Graph-theoretic iterative methods.} Spielman and Teng \cite{spielman-teng} and
successors \cite{kelner-stoc} build spectral sparsifiers as preconditioners, achieving
for any SDD system a provably near-linear $\mathcal{O}(m\,\mathrm{polylog}\,n\cdot\log(1/\varepsilon))$ bound. The hidden constants are large; the fast practical realization is the
approximate-Cholesky (approxChol) solver \cite{kelner-stoc,laplacians}, recently made
substantially more robust by Gao, Kyng, and Spielman (GKS) \cite{gks2023}, and paralleled by randomized-Cholesky direct
solvers such as RCHOL \cite{rchol}. We benchmark against both approxChol variants as
our strongest near-linear baselines (\S\ref{sec:compete}).

\emph{Algebraic multigrid.} AMG \cite{brandt-amg,ruge-stuben} builds a hierarchy of coarser
operators from matrix entries alone. Classical AMG selects a coarse node subset; aggregation
AMG \cite{vanek-mandel-brezina,notay} groups fine nodes into aggregates. AMG excels on
discretized PDEs but its complexity is hard to control on general graphs: coarse operators
densify on scale-free hubs, and higher-order interpolation makes this worse. Bootstrap AMG
\cite{brandt-bootstrap} and algebraic distance \cite{ron-safro-brandt} learn the coarsening from
relaxed test vectors, at higher setup cost; a recent direction learns the prolongation itself with
graph neural networks \cite{learned-amg}. LAMG \cite{lamg}, by contrast, is parameter-free and
training-free.

Generic AMG for graph Laplacians has grown into an active, rapidly expanding line:
Napov and Notay's degree-aware rooted aggregation with quality control \cite{napov-notay} (the
closest algorithmic relative; discussed in \S\ref{sec:draqc}), matching-based multilevel
preconditioners \cite{brannick-match}, effective-resistance/diffusion-distance coarsening with
least-squares interpolation from relaxed vectors \cite{lee2024}, Combinatorial Multigrid (CMG)
\cite{cmg}, and a growing body of surveys and applications spanning coarsening for machine
learning \cite{chen-saad-coarsening}, renormalization of complex networks
\cite{laplacian-coarsegrain}, and image-analysis solvers \cite{he-image-amg}. General-purpose
packages (PyAMG \cite{pyamg}, hypre/BoomerAMG \cite{boomeramg}) are tuned for PDE
discretizations rather than the irregular degree distributions of graph Laplacians. In the
broadest recent comparison \cite{gks2023}, CMG, BoomerAMG, and PETSc all failed to converge
across the chimera and SDDM families, leaving the approximate-elimination solver AC (the
\texttt{approxchol\_lap2} variant) the only blackbox method robust across all of them; the
original MATLAB LAMG was reported non-convergent there, whereas our independent re-derivation
converges on those families (\S\ref{sec:families})---a failure specific to that evaluation, not
the algorithm. We thus consider approxChol as our sole robust near-linear baseline, with CHOLMOD
for the direct-method crossover. We also compare with the other AMG solvers for the specific
graph classes they converge on.

\subsection{Contribution}\label{sec:contrib}
LAMG \cite{lamg} is an accelerated caliber-1 aggregation AMG that attains $\Om$ on general
graphs---more precisely $\mathcal{O}(m\log(1/\varepsilon))$ for an $\varepsilon$-accurate
solution---by regulating the coarse-operator \emph{energy} rather than aggregate size. We present
\textbf{LAMG+}: an open re-derivation of the algorithm\footnote{Open-source Julia implementation,
including all scripts to reproduce the experiments in this paper: \url{https://github.com/orenlivne/lamgplus}.}
(validated line-by-line against the
original MATLAB reference) that we use to establish three things.
\begin{itemize}[leftmargin=1.4em,itemsep=1pt]
\item \emph{A fair, same-language adjudication against the modern near-linear champion}
(\S\ref{sec:compete}). LAMG (2012) predates and was never compared with GKS's approxi\-mate-Cholesky
line \cite{laplacians,gks2023}, now the dominant practical near-linear solver.
Run in the same Julia process, the two are \emph{complementary peers}: approxChol is
faster on social/citation graphs, LAMG+ on FE/structural matrices (the fastest robust solver there,
$2.2\times$ over the robust AC variant on large graphs). The split is algorithmic
(coarsening vs.\ fill), not a compiler artifact.
\item \emph{Linear scaling at scale} (\S\ref{sec:scaling}). LAMG+ is empirically $\Om$ over
$1{,}711$ SuiteSparse graphs---$100\%$ converged, median $4$ cycles, log--log slope $1.01$ up to $2.4\!\times\!10^8$ nonzeros (twice larger than LAMG's original demonstration) with no tuning.
\item \emph{Algebraic diagnosis and resolution of LAMG's outlier tail} (\S\S\ref{sec:choices},\,\ref{sec:families}).
Where the original MATLAB LAMG was reported by \cite{gks2023} as non-convergent across the benchmark's
test families, we run the \emph{unmodified} LAMG~2.2.1 on those families (regenerated from GKS's released code) with the benchmark's exact solver call and find it converges on every one, as does our
re-derivation to the $10^{-8}$ tolerance (Table~\ref{tab:lamgorig})---placing the reported failure
in the evaluation, not the algorithm.
Additionally, by deploying a Local Fourier Analysis (Appendix~\ref{sec:aniso}), we prove that caliber-1 interpolation forces a strict energy leak on grid-aligned anisotropy. We resolve this with two lean, \emph{local} refinements---a strength-of-connection veto and a selective caliber-2
interpolation---turning LAMG's stalled
grid-aligned-anisotropy rate (Asymptotic Convergence Factor $(\acf)\approx0.99$) into $\approx0.11$ at $1\%$ median
overhead, each acting only on the locally anisotropic nodes that need it.
\item \emph{An ablation study of LAMG's load-bearing design choices} (\S\S\ref{sec:alg},\,\ref{sec:ablation}).
Several of LAMG's setup choices---no aggregate-size cap, orphan nodes left as singleton seeds, a
directional energy guard---are individually load-bearing yet undocumented as essential; we
isolate them and quantify, by single-choice ablation, that replacing any one by a plausible correlate
silently destroys convergence on a whole graph class (web ACF $0.002\!\to\!0.97$, for one). Here
reproduction \emph{is} discovery: knowing which heuristics carry the algorithm, and by how much, is
what makes it transferable beyond the original implementation.
\end{itemize}

\section{The LAMG algorithm}\label{sec:alg}
We summarize LAMG \cite{lamg} in enough detail to make this paper self-contained; the
refinements of \S\ref{sec:choices} presuppose it.

\paragraph{Multigrid in one paragraph} A relaxation such as Gauss--Seidel (GS) applied to
$L\phi=b$ removes the high-residual (oscillatory) error in a few sweeps but stalls on the
\emph{algebraically smooth} error---the slow, low-energy modes. Multigrid corrects those
on a coarser graph: it represents the fine smooth error $e\approx Pe^c$ by interpolation
$P\in\R^{n\times n_c}$ from a coarse vector $e^c$, solves the Galerkin coarse system
$L^c e^c=P^{\top}(b-L\tilde\phi)$ with $L^c=P^{\top}LP$ and $\tilde\phi$ the approximation to $\phi$ post-relaxation, and corrects
$\tilde\phi\leftarrow\tilde\phi+Pe^c$. Applied recursively over $\ell=1,\dots,L$ levels
this is a multigrid cycle; the work is geometric in the coarsening ratio and therefore typically linear in $m$.
LAMG uses two coarsening operators---elimination and aggregation---alternated until the
coarsest level is small.

\paragraph{Relaxation} LAMG employs Gauss--Seidel,
$\phi_u\leftarrow\big(b_u-\sum_{v\neq u}L_{uv}\phi_v\big)/L_{uu}$, an effective smoother
for SPD systems that needs no damping parameter.

\paragraph{Low-degree elimination} Let $F\subseteq V$ be a maximal independent set of
nodes of degree $\le d_{\max}=4$. Because $F$ is independent, $L_{FF}$ is diagonal and the
Schur complement onto $C=V\setminus F$ is \emph{exact}:
\begin{equation}\label{eq:elim}
L^c=L_{CC}-L_{CF}L_{FF}^{-1}L_{FC},\qquad
\phi_F=P_F\phi_C+L_{FF}^{-1}b_F,\quad P_F=-L_{FF}^{-1}L_{FC}.
\end{equation}
Eliminating an $F$-node of degree $\le 3$ does not increase $m$; the rare degree-$4$ node
may add up to two edges. Elimination removes the effectively-one-dimensional part of the graph at
zero approximation cost (up to fill), shrinking $n$ before each aggregation step.

\paragraph{Caliber-1 aggregation} Each remaining node is assigned to exactly one
aggregate, giving a piecewise-constant (``caliber-1'') interpolation with $P_{u,A(u)}=1$.
This keeps $L^c=P^{\top}LP$ \emph{as sparse as the graph}---it cannot manufacture the dense
coarse rows that sink classical AMG on hubs (hence the ``lean'' in LAMG). Aggregates are formed from $K$ \emph{test
vectors} $X=(x^{(1)},\dots,x^{(K)})$, each being the result of $\nu$ GS sweeps on $Lx=0$ from a random start:
these are samples of the algebraically smooth subspace the coarse grid must capture, with
no geometry assumed. (LAMG uses $K=4$ at the finest level, raises $K$ by one per coarser
level up to $10$, and smooths each by $\nu=3$ sweeps \cite[\S4.1]{lamg}.) The proximity of an edge $(u,v)$ is the \emph{affinity}, the squared
cosine of the two test-vector profiles \cite[Eq.~3.3]{lamg},
\begin{equation}\label{eq:aff}
c_{uv}=\frac{\big(\sum_k x^{(k)}_u x^{(k)}_v\big)^2}
{\big(\sum_k (x^{(k)}_u)^2\big)\big(\sum_k (x^{(k)}_v)^2\big)}\in[0,1].
\end{equation}
A node $u$ joins the seed $s$ of maximal affinity that is \emph{admissible} under the
\emph{energy-ratio guard}. With the nodal energy
$E_u(x;y)=\tfrac12 L_{uu}y^2-B_u(x)\,y+C_u(x)$, $B_u=\sum_v w_{uv}x_v$,
$C_u=\tfrac12\sum_v w_{uv}x_v^2$,
\begin{equation}\label{eq:guard}
q_{u\leftarrow s}=\max_{1\le k\le K}
\frac{E_u\!\big(x^{(k)};x^{(k)}_s\big)}{\min_y E_u\!\big(x^{(k)};y\big)}\le Q,\qquad Q=2.5 .
\end{equation}
The numerator is $u$'s energy when caliber-1 \emph{forces} $u=x_s$; the denominator is its
relaxed (GS-optimal) energy $\min_y E_u=E_u(x;B_u/L_{uu})$. Their ratio is the factor by
which lumping $u$ into $s$'s aggregate inflates the local energy; capping it bounds the
coarse-operator energy aggregate by aggregate, hence the two-level factor
$\rho\lesssim 1-1/Q$ \cite[\S3.4]{lamg}. Crucially the guard regulates \emph{energy, not
size}---there is no hard aggregate-size cap: the same $Q$ admits a $100$-node hub
aggregate (neighbors agree on the slow modes, $q\approx1$) yet forbids a $3$-node mesh
aggregate (neighbors disagree, $q$ shoots up). Two companion conventions complete the rule
and are equally load-bearing: a node that finds no admissible seed stays a singleton aggregate;
and the guard is \emph{directional}---\eqref{eq:guard} scores $u$ forced onto seed $s$, $q_{u\leftarrow s}$, and is not symmetrized in $(u,s)$. These
are LAMG's choices \cite[\S3.4]{lamg}; substituting a plausible correlate for any of them (a
size cap, forced absorption, a symmetric guard) silently over- or under-coarsens some graph
class, so \S\ref{sec:choices} discusses where LAMG+ deviates from LAMG. The guard
dominates setup cost and is evaluated in $\Om$ via three bulk-precomputed moments per node
(App.~\ref{sec:guardeval}).

\paragraph{Cycle and acceleration} The solve runs a $\gamma$-cycle whose per-level index
LAMG fixes by a work rule \cite[Eq.~4.1]{lamg}: elimination levels use $\gamma_\ell=1$, and an
aggregation level uses
\begin{equation}\label{eq:gamma}
\gamma_\ell=\begin{cases}1.5, & |E_\ell|>0.1\,|E|\ \ \text{(finer levels)},\\[2pt]
\min\{2,\;0.7/\tau_\ell\}, & \text{otherwise (coarser levels)},\end{cases}
\end{equation}
where $\tau_\ell=|E_{\ell+1}|/|E_\ell|$ is the level-$\ell$ coarsening ratio. The cap bounds
per-cycle work geometrically, guaranteeing $\Om$. The coarse-level branch raises $\gamma_\ell$
toward $2$ where coarsening is aggressive but lowers it toward $1$ where coarsening is weak
($\tau_\ell\to1$); LAMG+ keeps the same rule with a slightly looser work cap.
After each cycle the finest iterate is improved by a min-residual
(Krylov) recombination over a short history of $\kappa$ saved iterates \cite[\S7.8.2]{guide}:
$\alpha=\arg\min_\alpha\|r-LE\alpha\|$ with $r = b - L \tilde\phi$, $E=[\phi_1-\tilde\phi,\dots,\phi_{\kappa}-\tilde\phi]$, then
$\tilde\phi\leftarrow \tilde\phi+E\alpha$. A level's history accumulates across the cycle's $\gamma$
revisits; the finest persists across cycles.

\section{Two refinements beyond LAMG}\label{sec:choices}
Reproducing LAMG (\S\ref{sec:alg})---no size cap, singleton nodes, a directional guard, the
cycle-index calibration---is a prerequisite, not a contribution. LAMG+ departs
from LAMG (Table~\ref{tab:diff}) by (1)~avoiding aggregation across a weak edge, via a
strength-of-connection veto (\S\ref{sec:soc}), and (2)~judiciously raising the interpolation caliber
from $1$ to $2$ (\S\ref{sec:cal2}); both target (local) grid-aligned anisotropy---the bounded-but-slow
outlier tail LAMG reported \cite[\S\S5.2,\,6.1]{lamg}. Together they turn a stalled
$\acf\approx0.99$ into $\approx0.11$ on the worst anisotropic grids while leaving every other class
untouched. Our contribution is their lean, parameter-free realization and
at-scale evidence (\S\ref{sec:results}). Both refinements are \emph{local}---the veto a per-edge
test, caliber-2 a per-node test, each firing only where the matrix is locally anisotropic and inert
($\approx\!0\%$ of nodes) elsewhere---so refinement follows where the solver is locally slow, never a
global average: a graph that is mostly isotropic with a small anisotropic region defeats any global
trigger yet is handled correctly, node by node (\S\ref{sec:calib}). Because the two act at different
stages (which nodes aggregate, how they interpolate), they are independent and individually
ablatable (\S\ref{sec:ablation}).

\begin{table}[h!]\centering
\caption{The two places LAMG+ departs from LAMG, both targeting grid-aligned anisotropy. Reproduced setup choices (no size cap, singletons, directional guard,
non-decaying $\gamma$, fixed $K{=}4$) are discussed in \S\S\ref{sec:alg},\,\ref{sec:calib}.}\label{tab:diff}
\begin{tabular}{>{\raggedright\arraybackslash}p{2.1cm}
                >{\raggedright\arraybackslash}p{2.4cm}
                >{\raggedright\arraybackslash}p{3.5cm}
                >{\raggedright\arraybackslash}p{3.1cm}}
\toprule
Refinement & LAMG & LAMG+ & Effect\\
\midrule
Aggregation & affinity $+$ energy guard only & $+$ strength-of-connection veto: no merge across
a weak edge (eq.~\eqref{eq:soc}) & removes the rare affinity-inversion stall,
$\acf\ 0.99\!\to\!0.46$\\
Interpolation caliber & strictly caliber-1 & $+$ selective caliber-2 on locally-$1$-D nodes &
lowers the $q\!\approx\!2$ ceiling, $\acf\ 0.46\!\to\!0.11$\\
\bottomrule
\end{tabular}
\end{table}

\subsection{A strength-of-connection veto}\label{sec:soc}
Affinity \eqref{eq:aff} ranks a node's neighbors by how well their relaxed test-vector
profiles align, and aggregation merges a node with its highest-affinity admissible neighbor.
On a strongly grid-aligned anisotropic grid the test vectors are nearly flat along the weak
direction, so the $K$-sample affinity of a weak (weight-$\epsilon$) edge is computed from
almost-identical profiles and can, on an unlucky sample, exceed that of the genuinely
strong (weight-$1$) edge---we measure $0.97$ versus $0.89$ at one interior node of a $256^2$,
$\epsilon=10^{-4}$ grid. Neither the affinity nor the energy guard $Q$ catches this: both are
computed from the \emph{same} test vectors and are fooled identically. Fortunately, the matrix is
not in this case---the conductance ratio there is $10^4$. A single such inversion creates one
cross-weak-direction merge, a localised residual hot-spot (peak-to-mean $>10^3$), and a stalled
rate ($\acf\approx0.99$) that no amount of $\gamma$ or recombination repairs, and that
caliber-2 alone cannot reach past either.

The veto is the classical Ruge--St\"uben strength-of-connection criterion \cite{ruge-stuben}
applied to aggregation: never aggregate across an edge whose weight is below a fraction
$\tau_{\mathrm{soc}}$ of the node's strongest incident weight,
\begin{equation}\label{eq:soc}
\text{skip edge }(u,v)\ \text{in aggregation if}\quad |w_{uv}|<\tau_{\mathrm{soc}}\,\max_{t\sim u}|w_{ut}|,
\qquad \tau_{\mathrm{soc}}=0.05 .
\end{equation}
It costs $\Oone$ per edge, changes no asymptotics, and is inert wherever a node's incident
weights are comparable. With it, the $256^2$ and $384^2$ grids drop from $\acf\approx0.99$ to
$0.11$ across all sizes, seeds, and $\epsilon\le10^{-1}$ (worst $0.117$ over $24$ cases); on a
$42$-graph random SuiteSparse sample and the social/FE giants to $1.5\times10^7$ edges it is
neutral (worst ACF ratio $1.22$, no regression). The defect is rare---it needs a near-degenerate
weak direction \emph{and} an unlucky sample---which is why stock LAMG mostly avoids it and
reports anisotropy as bounded-but-slow rather than stalled; the veto removes the residual tail
and, crucially, is what lets the caliber-2 fit of \S\ref{sec:cal2} act on a clean aggregation.

\subsection{Selective caliber-2 on one-dimensional chains}\label{sec:cal2}
The veto fixes \emph{which} nodes aggregate; caliber-2 changes \emph{how} they
interpolate---the place LAMG's strict caliber-1 design leaves a bounded but non-textbook
rate. On grid-aligned anisotropy the affinity correctly semicoarsens along the strong direction
(the aggregates are well placed), but the caliber-1 piecewise-constant interpolation across each
$2{:}1$ aggregate suffers a structural order deficit. As proven by the Local Fourier Analysis
in Appendix~\ref{sec:aniso}, caliber-1 commits exactly half of every interpolant's energy to
the unresolved high-frequency harmonic ($\beta_{c1}=1/2$), inflating the energy ratio to
$q\approx2$ and pinning the two-level factor at $\rho\approx0.5$. By raising the interpolation
to caliber-2 exclusively where a fine node's strong edges (i.e., those not vetoed
by~\eqref{eq:soc}) reach exactly two coarse aggregates, the symbol's dependence on the smooth
mode shifts from
$O(\theta_x)$ to $O(\theta_x^2)$, mathematically closing the energy leak ($\beta_{c2}\to1$).
LAMG+ implements this by adding a second, caliber-2 parent with a single fitted weight $w$
(App.~\ref{sec:aniso} derives it). Measured per level with an exact coarse solve (so the
$\gamma$ schedule cannot mask it), the two-level factor on the $256^2$, $\epsilon=10^{-4}$ grid
falls from $0.46$ (caliber-1, veto on) to $0.28$ at every aggregation level; with $\gamma$ and
recombination the full-cycle rate reaches $\approx0.11$. The gate is \emph{self-targeting}:
$\approx\!0\%$ of nodes upgrade on isotropic graphs, so over a $1{,}267$-graph SuiteSparse
sample it rescues the anisotropic minority (\texttt{NACA0015}, \texttt{epb3}, \texttt{darcy}
from $100$ cycles to $5$--$8$; the \texttt{thermomech}/\texttt{thermal} stiffness matrices from
non-convergence to $\approx\!13$) while leaving the easy majority's solve essentially unchanged,
at $\approx\!1\%$ median wall-clock. The two refinements are complementary:
caliber-2 without the veto still stalls at $0.99$ on the worst grids, the veto without caliber-2
reaches only the caliber-1 ceiling $\approx0.46$; together they give $0.11$ (\S\ref{sec:ablation}).

\subsection{A single, local configuration}\label{sec:calib}
LAMG+ runs one configuration on every graph: a cheap $K{=}4$ affinity with both refinements always
on. The interpolation caliber is raised from $1$ to $2$ wherever a node is locally one-dimensional
(\S\ref{sec:cal2}) and left at caliber-1 everywhere else, so the extra cost falls only on the
anisotropic nodes. Four test vectors are an ample sample for the single least-squares weight each
upgraded node fits, so we simply fix $K{=}4$ rather than accumulate test vectors with depth as the
original LAMG does (raising $K$ helps reduce the probability of a spurious affinity, but would significantly pay off only together with more relaxation sweeps $\nu$ to increase the interpolation caliber beyond $2$, which we avoid in this lean AMG framework).

Because the caliber decision is taken per node, it needs no global control, making it robust. A graph that is mostly isotropic with a small anisotropic region is handled correctly:
the upgrade fires exactly on that region. In a $64^2$
anisotropic patch embedded in a $256^2$ or $512^2$ isotropic grid, only $6\%$ and $1.5\%$ of the
nodes trigger, respectively (Table~\ref{tab:adaptive}). The caliber-1 $\acf\approx0.3$ improves to a textbook multigrid efficiency $\acf\approx0.02$ for caliber-2.

\begin{table}[t]\centering
\caption{Caliber-1 interpolation vs. caliber-2 ($\acf$, with cycles to $\varepsilon=10^{-8}$ in parentheses). ``aniso.\ nodes'' is the fraction of locally
one-dimensional nodes. Non-anisotropic graphs are unaffected (\S\ref{sec:ablation}).}\label{tab:adaptive}
\begin{tabular}{lccc}
\toprule
graph & aniso.\ nodes & caliber-1 & caliber-2 (default)\\
\midrule
aniso $256^2$, $\epsilon{=}10^{-4}$        & $99\%$  & 0.56 (13) & \textbf{0.16 (7)}\\
aniso $384^2$, $\epsilon{=}10^{-4}$        & $99\%$  & 0.83 (33) & \textbf{0.11 (6)}\\
isotropic $256^2$ $+$ $64^2$ aniso.\ patch & $6\%$   & 0.31 (8)  & \textbf{0.02 (5)}\\
isotropic $512^2$ $+$ $64^2$ aniso.\ patch & $1.5\%$ & 0.24 (7)  & \textbf{0.02 (5)}\\
isotropic grid $256^2$                     & $0\%$   & 0.016 (5) & 0.015 (5)\\
\texttt{bodyy5} (FE stiffness)             & $0\%$   & 0.020 (5) & 0.019 (5)\\
\bottomrule
\end{tabular}
\end{table}

\section{Numerical results}\label{sec:results}
LAMG+ is implemented in Julia, and all reported times are wall-clock, measured on an Apple M5~Pro
(18-core, 48~GB RAM); the competitor solvers of \S\ref{sec:compete} run in the same Julia process
under identical compilation, so the wall-clock comparisons are language-fair.
Throughout this section, every solve uses the right-hand side $b=L\,x_{\mathrm{true}}$ with
$x_{\mathrm{true}}$ a standard-normal random vector projected to zero mean, and runs to
$\varepsilon=\|b-L\,x\|/\|b\|\le10^{-8}$ unless a different $\varepsilon$ is stated.

\paragraph{Test corpus} All experiments draw from one fixed corpus, the union of two sources:
(i)~the complete $1{,}711$-graph SuiteSparse collection used as the LAMG test set
($100\le n\le1.2\!\times\!10^7$), which drives the scaling study (\S\ref{sec:scaling}) and, restricted
to its $m>10^6$ members, the head-to-head competition (\S\ref{sec:compete});
(ii)~instances from the AC robustness benchmark \cite{gks2023}, generated by its
own generators---all $30$ \texttt{chimera} and weighted-\texttt{chimera} graphs and all $9$
uniform/anisotropic/high-contrast Poisson grids---plus our synthesized Kyng--Sachdeva classes
(Erd\H{o}s--R\'enyi, Barab\'asi--Albert, stochastic-block, $k$-regular, power-law-cluster), used for
robustness (\S\ref{sec:families}). Every solver runs on the same exported $(L,b)$; all per-graph setup and solve
times are recorded and released with the code \cite{lamgplus-code} so the tables can be regenerated. The largest graph
solved end-to-end is \texttt{com-Orkut} ($1.2\!\times\!10^8$ edges, $2.4\!\times\!10^8$ nonzeros,
$2$ cycles).

\subsection{Convergence against the reference}\label{sec:acf}
Table~\ref{tab:acf} reports the asymptotic convergence factor (ACF, the geometric mean of
the last five per-cycle residual reductions) on identical exported systems $(L,b)$. LAMG+
matches or beats the MATLAB reference on every class.

\begin{table}[t]\centering\small
\caption{ACF on identical $(L,b)$; ``cyc'' is cycles to $\varepsilon=10^{-8}$.
``LAMG+'' and ``ref.'' are the LAMG+ and MATLAB-reference asymptotic convergence factors.
Entries $<10^{-4}$ solve in $1$--$2$ cycles---below the resolution at which an asymptotic
factor is meaningful. Class abbreviations: AS = autonomous-systems, P2P = peer-to-peer.}\label{tab:acf}
\begin{tabular}{llrrccc}
\toprule
Graph & Class & $n$ & $m$ & LAMG+ & cyc & ref.\\
\midrule
grid128        & structured  & 16\,384 & 32\,512  & 0.022      & 5 & 0.21\\
bodyy5         & FE stiff.   & 18\,589 & 55\,132  & 0.019      & 5 & 0.27\\
wb-cs-stanford & web         & 7\,008  & 17\,054  & 0.0009     & 3 & 0.06\\
ca-GrQc        & collab.     & 4\,158  & 13\,422  & 0.0008     & 3 & 0.07\\
ca-HepTh       & collab.     & 8\,638  & 24\,806  & 0.0040     & 4 & 0.08\\
Oregon-1       & AS          & 11\,174 & 23\,409  & $<10^{-4}$ & 2 & 0.005\\
as-caida       & AS          & 26\,475 & 53\,381  & $<10^{-4}$ & 2 & 0.008\\
p2p-Gnutella08 & P2P         & 6\,262  & 9\,719   & 0.0010     & 3 & 0.03\\
email-Enron    & social      & 33\,696 & 180\,811 & 0.0011     & 3 & 0.10\\
ak2010         & road        & 42\,381 & 102\,091 & 0.0066     & 4 & 0.12\\
\bottomrule
\end{tabular}
\end{table}

\subsection{Linear scaling}\label{sec:scaling}
\begin{figure}[t]\centering\includegraphics[width=\linewidth]{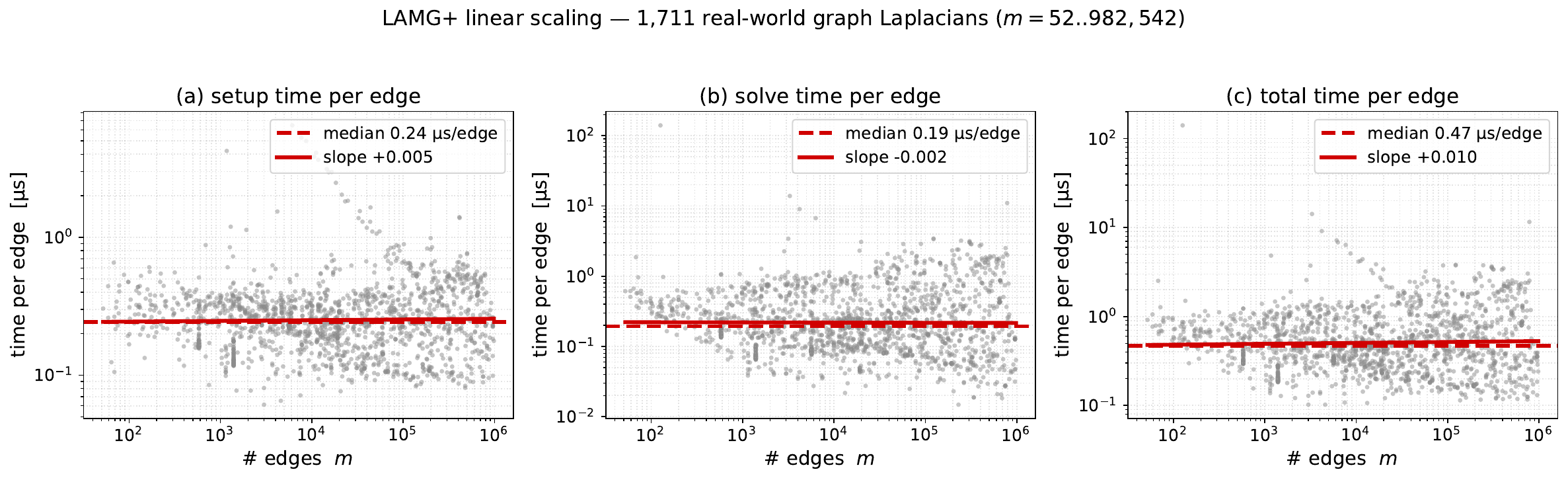}
\caption{LAMG+ per-edge time vs.\ number of edges $m$ (log--log) over $1{,}711$ SuiteSparse graph
Laplacians: (a)~setup, (b)~solve, (c)~total. Grey points are individual graphs; the red dashed line is
the median, the red solid line the log--log regression fit (slope shown per panel). A near-zero slope
means per-edge cost is independent of $m$ --- the empirical $\Om$ statement; the regression of total
setup$+$solve time gives the exponent $m^{1.01}$.}\label{fig:scaling}\end{figure}

Figure~\ref{fig:scaling} plots LAMG+ per-edge setup and solve time versus $m$ for all $1{,}711$
SuiteSparse graphs ($100\le n\le1.2\!\times\!10^7$). Every one converges to $\varepsilon=10^{-8}$, in a median of $4$ cycles. The log--log slopes of per-edge time are $+0.005$ (setup)
and $-0.002$ (solve): per-edge cost is independent of $m$, supporting the empirical $\Om$ statement
(equivalently, total setup$+$solve time scales as $m^{1.01}$). The medians are
$t_{\text{setup}}/m\approx0.24~\mu$s/edge (IQR $[0.18,0.34]$) and
$t_{\text{solve}}/m\approx0.19~\mu$s/edge (IQR $[0.10,0.43]$). These low, flat constants reflect an
implementation tuned to the algorithm: a self-gated reverse-Cuthill--McKee reordering~\cite{cuthill-mckee} for cache
locality, single-pass Gustavson-fused Schur extraction (each coarse operator built in one column scan
with no matrix-multiply temporary), allocation-free elimination into preallocated buffers, and a
monomorphic, branch-free, residual-maintaining in-place Gauss--Seidel kernel. Each of these features does not change the result but lowers the hidden constant.
The largest per-edge costs are the strongly anisotropic grids; cf.~App.~\ref{sec:aniso}.

\subsection{Comparison: complementary strengths by class}\label{sec:compete}
We benchmark setup$+$solve to $\varepsilon=10^{-8}$ against both approximate-Cholesky
variants---the original, faster \texttt{approxchol\_lap} (``approxChol'') \cite{laplacians} and the
robust \texttt{approxchol\_lap2} (``AC'') of \cite{gks2023}, and note sparse Cholesky
(CHOLMOD) \cite{davis-cholmod} for context. All run in the same Julia process under identical
compilation, so any difference is algorithmic (coarsening vs.\ fill), not a compiler artifact; we
report time per nonzero ($\mu$s/nnz, the worst-case metric of \cite{gks2023}), each a warm timed run
(a warm-up call first excludes Julia's one-time JIT compilation), measured serially so the per-graph
rows and the aggregate below are contention-free and mutually consistent.

\begin{table}[t]\centering
\caption{Setup$+$solve time per nonzero ($\mu$s/nnz) to $\varepsilon=10^{-8}$; bold $=$
fastest. ``approxChol'' is the fast \texttt{approxchol\_lap}; ``AC'' the robust
\texttt{approxchol\_lap2}. Sparse Cholesky is omitted: it is fastest on the web/road graphs but
\emph{cannot factor} the FE/structural and dense-citation matrices ($>\!21$\,GB fill on
\texttt{flickr}). The bottom block summarizes LAMG+ against each competitor over \emph{all} $201$
graphs with $m>10^6$ ($\mathrm{nnz}(L)\le6\!\times\!10^7$): the \emph{geometric-mean speedup} is the
geometric mean of (competitor time)$/$(LAMG+ time), so $>1$ means LAMG+ is faster; the \emph{win
rate} is the fraction of the $201$ graphs on which LAMG+ is faster; each reported at $\varepsilon=10^{-8}$\,/\,$10^{-4}$. All times are serial (contention-free).}\label{tab:compete}
\begin{tabular}{llrrrr}
\toprule
Graph & family & edges & LAMG+ & approxChol & AC\\
\midrule
web-Stanford     & web        & 1.0\,M  & 0.26 & \textbf{0.19} & 0.30\\
web-Google       & web        & 2.6\,M  & 0.48 & \textbf{0.25} & 0.43\\
roadNet-PA       & road       & 1.5\,M  & 0.72 & \textbf{0.26} & 0.34\\
flickr           & social     & 4.9\,M  & 0.55 & \textbf{0.26} & 0.62\\
coPapersDBLP     & citation   & 15.2\,M & 0.31 & \textbf{0.23} & 0.46\\
\midrule
bmwcra\_1        & FE         & 5.2\,M  & \textbf{0.08} & 0.28 & 0.57\\
crankseg\_1      & structural & 5.3\,M  & \textbf{0.08} & 0.19 & 0.65\\
troll            & structural & 5.9\,M  & \textbf{0.10} & 0.24 & 0.44\\
pwtk             & structural & 5.7\,M  & \textbf{0.11} & 0.21 & 0.40\\
\midrule
\multicolumn{6}{l}{\emph{all} $m>10^6$ graphs ($201$), LAMG+ vs.\ each ($10^{-8}$\,/\,$10^{-4}$):}\\
\multicolumn{4}{l}{\quad geometric-mean speedup} & $0.99$\,/\,$1.20\times$ & $2.20$\,/\,$2.91\times$\\
\multicolumn{4}{l}{\quad win rate}               & $56$\,/\,$64\%$        & $76$\,/\,$84\%$\\
\bottomrule
\end{tabular}
\end{table}

The finding is a clean \emph{class split}, not a single winner (Table~\ref{tab:compete}), but at root it is one
mechanism---a \emph{fill contest}. These graphs all converge in only $3$--$4$ cycles, so neither
solver is paced by iteration count; the winner is whichever builds the leaner hierarchy,
because the fill a solver creates sets both its setup cost and the cost of every cycle. Aggregation
and sampled elimination have opposite fill profiles. On low-diameter social, web, and citation
graphs and on the planar road network, approxChol wins: their low mean degree keeps its sampled
elimination cheap, while the absence of clean separators inflates LAMG+'s coarse operators---and with
them the cost of every cycle. On finite-element and structural matrices the profiles reverse---the high per-node
coupling makes elimination fill in, while LAMG+'s aggregation collapses those dense neighborhoods at
$\Oone$ cost on a cleanly separating mesh---and so does the order: LAMG+ is the fastest of the three (geometric-mean
$1.76\times$ over approxChol across the $34$ FE/structural graphs, winning $91\%$ of them;
\texttt{bmwcra\_1}, \texttt{crankseg\_1}, \texttt{troll}, \texttt{pwtk} run $1.9$--$3.5\times$ over
it), the leanest build of the three, and the fixed cheap $K{=}4$ calibration (\S\ref{sec:calib})
pays off. Aggregated over all
$201$ large graphs ($m>10^6$), the two are \emph{near-peers}: LAMG+ wins $56\%$ of head-to-heads
against approxChol (geometric-mean $0.99\times$---approxChol's social/web
wins are by similar factors to LAMG+'s FE/structural wins), and $76\%$ against the robust AC, where
it is $2.2\times$ faster in geometric mean. (At the lower $10^{-4}$ accuracy many applications need, the balance shifts toward LAMG+; cf.~\S\ref{sec:tol4}.) Sparse Cholesky is fastest where good separators exist (web crawls, road)
but cannot factor the largest poorly-separable matrices at all.

Table~\ref{tab:classcomp} widens the comparison to all seven solvers across 13 graph classes,
including all of the AC benchmark's synthetic families. The result is a clean
separation by design intent. On the structured classes they were built for---FE/structural
and grids---the four algebraic-multigrid solvers (LAMG+, pyAMG, BoomerAMG, PETSc GAMG) are all fast
and beat sparse Cholesky. On the low-diameter classes (social, citation, web) BoomerAMG~\cite{boomeramg} and PETSc GAMG slow by up to an order of magnitude, pyAMG~\cite{pyamg}
stops converging, and approxChol is fastest. CMG~\cite{cmg} does not converge on the grids,
on SPE, or on many roads. BoomerAMG converges everywhere but is catastrophic ($7\,\mu$s/nnz) on
irregular optimization graphs; pyAMG, CMG, and PETSc GAMG additionally fail to converge on $11$, $15$, and $6$ of
the instances, respectively. Across all classes only three solvers are both convergent
everywhere and never catastrophic (worst case under $1.4\,\mu$s/nnz): \textbf{LAMG+}
($\le1.35$), approxChol ($\le0.70$), and AC ($\le1.20$)---precisely the robust near-linear
solvers this paper compares.

\begin{table}[t]\centering\footnotesize\setlength{\tabcolsep}{3.4pt}
\caption{Multi-solver comparison by graph class: mean\,/\,worst-case setup$+$solve time per nonzero
($\mu$s/nnz) to $\varepsilon=10^{-8}$, over a class-stratified corpus (real SuiteSparse/SNAP
classes and \emph{all} of the approximate-Cholesky benchmark's synthetic families, including the
SPE10 reservoir). Bold $=$ fastest mean; ``DNF'' $=$ all instances failed to converge.
A per-graph budget of $50\times$ the fastest solver (capped at $90$\,s) bounds the slow cases.
apxChol $=$ \texttt{approxchol\_lap}; PETSc $=$ PETSc GAMG.}\label{tab:classcomp}
\begin{tabular}{lrrrrrrr}
\toprule
class & LAMG+ & apxChol & AC & Boomer & pyAMG & CMG & PETSc\\
\midrule
FE/structural   & 0.12/0.15 & 0.28/0.54 & 0.53/0.73 & 0.18/0.35 & \textbf{0.10}/0.13 & 0.18/0.21$^{2}$ & 0.12/0.20\\
mesh/grid       & 0.84/1.24 & 0.24/0.32 & 0.55/0.74 & \textbf{0.18}/0.27 & 0.38/0.54 & 1.46/1.46$^{2}$ & 0.23/0.33\\
social/citation & 0.47/0.69 & \textbf{0.25}/0.37 & 0.71/0.94 & 1.53/2.50 & 4.60/7.81$^{3}$ & 0.50/0.73 & 3.10/3.14$^{4}$\\
web             & 0.37/0.41 & \textbf{0.29}/0.38 & 0.39/0.51 & 1.43/2.59 & 5.45/5.45$^{3}$ & 1.27/2.05$^{1}$ & 2.76/3.59$^{2}$\\
road            & 0.68/0.85 & 0.39/0.50 & 0.53/0.65 & \textbf{0.35}/0.39 & 1.38/1.76 & 1.24/1.24$^{3}$ & 0.57/0.75\\
optimization    & 0.15/0.18 & 0.26/0.36 & 0.88/1.03 & 4.07/7.02 & \textbf{0.13}/0.13 & 0.23/0.33 & 2.04/3.47\\
chimera         & 0.40/0.91 & \textbf{0.14}/0.19 & 0.56/1.20 & 0.48/0.67 & 1.13/1.74$^{1}$ & 0.43/0.78 & 0.91/1.72\\
wtd-chimera     & 0.45/1.34 & \textbf{0.23}/0.70 & 0.39/0.56 & 0.46/0.63 & 2.04/3.25$^{1}$ & 0.33/0.76 & 0.96/1.86\\
sddm-chimera    & 0.30/0.51 & \textbf{0.10}/0.13 & 0.21/0.32 & 0.44/0.57 & 2.21/2.80$^{1}$ & 0.34/0.47 & 1.74/3.26\\
grid-3D         & 1.30/1.35 & 0.32/0.37 & 0.70/0.72 & \textbf{0.19}/0.19 & 0.37/0.41 & DNF & 0.25/0.25\\
aniso-grid      & 0.61/0.68 & 0.15/0.16 & 0.38/0.41 & \textbf{0.15}/0.15 & 2.38/2.38$^{2}$ & 1.22/1.48$^{1}$ & 1.01/1.49\\
star            & 0.05/0.05 & 0.03/0.03 & 0.24/0.40 & \textbf{0.02}/0.02 & 0.03/0.03 & 0.04/0.06 & 0.04/0.04\\
SPE/reservoir   & 1.09/1.22 & \textbf{0.29}/0.30 & 0.56/0.59 & 0.29/0.31 & 2.72/2.95 & DNF & 0.90/1.29\\
\bottomrule
\end{tabular}
\par\smallskip{\footnotesize\itshape Superscript $k$: $k$ instances in the class exceeded the
300-cycle limit; shown values are medians over converged instances only.}
\end{table}

\paragraph{Memory} The same split holds in space (Table~\ref{tab:memory}). We compare memory only
against the robust solvers (approxChol, AC): the non-robust competitors store their hierarchies in
C and are not directly measurable in bytes/nnz. On large FE/structural matrices LAMG+'s no-fill
multilevel hierarchy is the most frugal of the three ($21$--$30$ bytes/nnz vs.\
$33$--$55$ for the factorizations); on web/social graphs it is the most expensive. Across every graph the robust AC costs more memory than
approxChol, so AC's guarantee is paid in both time and space. The table
reports the \emph{persistent} hierarchy; the transient setup peak adds only a single level's
$K$ test vectors ($K{\cdot}n$ floats), discarded once that level is coarsened rather than accumulated
across levels, and the fixed $K{=}4$ default halves this peak versus LAMG's $K{=}8$
accumulation, which governs whether the largest graphs build within a fixed memory budget.

\begin{table}[t]\centering
\caption{Memory usage (bytes per nonzero of $L$) of the three robust solvers.
LAMG+ is the most frugal on FE/structural matrices, heavier on web/social; AC always
exceeds \texttt{approxchol\_lap}. Non-robust solvers (BoomerAMG, CMG, pyAMG, PETSc GAMG)
store hierarchies in C and are not included.}\label{tab:memory}
\begin{tabular}{lrrr}
\toprule
Graph & LAMG+ & approxChol & AC\\
\midrule
troll, pwtk (structural)   & \textbf{23} & 34 & 44\\
bmwcra\_1 (FE)             & \textbf{22} & 35 & 47\\
bone010 (FE, $2.3\!\times\!10^7$) & \textbf{30} & 43 & 55\\
flickr (social)            & 80 & \textbf{40} & 53\\
coPapersDBLP (citation)    & 49 & \textbf{32} & 32\\
web-Stanford (web)         & 62 & \textbf{32} & 32\\
\bottomrule
\end{tabular}
\end{table}

\subsection{Robustness across the corpus, and the closest algorithmic relative}\label{sec:draqc}
Over the full SuiteSparse corpus of \S\ref{sec:scaling}, built once and solved to $10^{-8}$, LAMG+
(caliber-2) is the only near-linear solver that converges on \emph{every} graph, and with the
tightest tail---at most $27$ cycles, versus $77$ PCG iterations for approxChol, whose single
non-convergence is a chemical-process matrix (Table~\ref{tab:draqc}).

\paragraph{Availability of competing implementations}
The graph-Laplacian solver most closely related to LAMG is the degree-aware rooted aggregation with
quality control (DRA-QC) of Napov and Notay \cite{napov-notay}: like LAMG it eliminates low-degree
vertices and coarsens by aggregation, but it controls aggregate quality through a provable two-grid
condition-number bound rather than an energy guard, and is reported $2$--$5\times$ faster on average
than the original LAMG on a real-world test set \cite{napov-notay}. AGMG \cite{notay}---the
general-purpose quality-control aggregation code---is distributed under a restrictive academic
licence, which we requested and were unable to obtain as a non-academic user; its DRA-QC variant is
not part of any public distribution and is covered by the same license terms. The other close relative---Lee's diffusion-distance/effective-resistance
affinity coarsening \cite{lee2024}, in spirit nearest to LAMG's relaxed-vector affinity---is likewise
not distributed in any public implementation. We therefore benchmark
quantitatively only against the publicly available approxChol \cite{laplacians,gks2023}---the
dominant practical near-linear solver.

\begin{table}[t]\centering
\caption{Robustness over the full SuiteSparse corpus (build the hierarchy once, solve to
$10^{-8}$): fraction of graphs converged and the per-graph iteration distribution. LAMG+
caliber-2 is the only solver convergent on every graph and has the tightest tail. ``it''~=~cycles
(LAMG+) or PCG iterations (approxChol).}\label{tab:draqc}
\begin{tabular}{lrrrrr}
\toprule
Solver & converged & median it & p99 it & max it & \#fail\\
\midrule
LAMG+ (caliber-2)         & \textbf{100.0\%} & 4  & 9  & \textbf{27}  & 0\\
approxChol \cite{gks2023} & 99.9\%           & 15 & 37 & 77           & 1\\
\bottomrule
\end{tabular}
\end{table}

\subsection{An a-priori selection rule}\label{sec:predictor}
\begin{figure}[t]\centering
\includegraphics[width=0.62\linewidth]{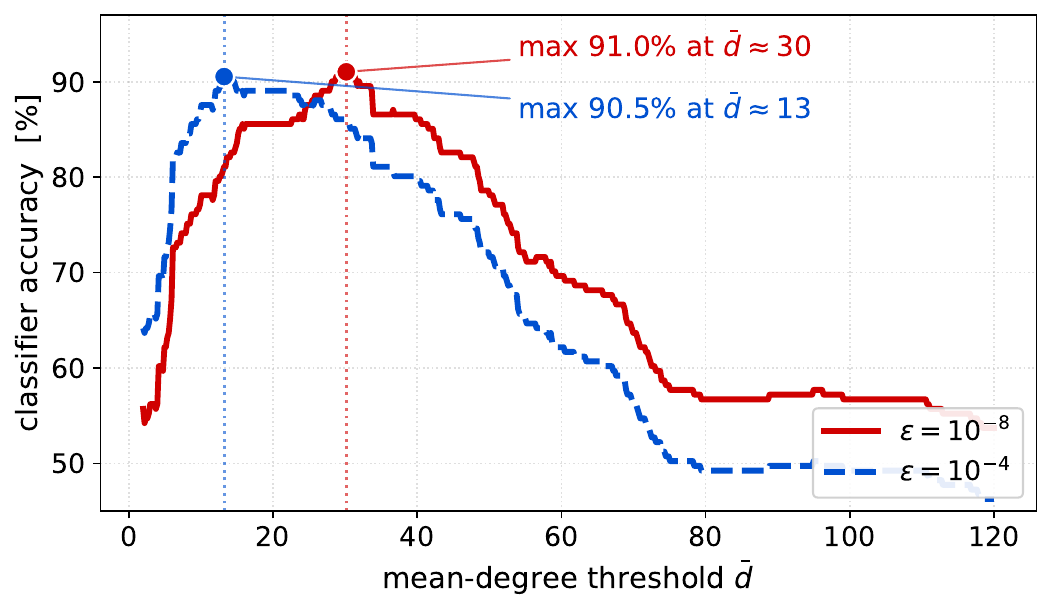}
\caption{Degree-threshold predictor of a LAMG+ win: accuracy of the rule ``LAMG+ wins iff
$\bar d>t$'' over the $201$ large ($m>10^6$) graphs as a function of the threshold $t$, at both
solve accuracies. At $\varepsilon=10^{-8}$ (red, solid) accuracy peaks at $91.0\%$ near
$\bar d\approx30$; at the looser $\varepsilon=10^{-4}$ (blue, dashed) LAMG+'s higher win rate
($64\%$ vs.\ $56\%$; \S\ref{sec:tol4}) shifts the optimal threshold left, peaking at $90.5\%$ near
$\bar d\approx13$. Both peaks are marked. The near-balanced $\varepsilon=10^{-8}$ classes ($112$
LAMG+ vs.\ $89$ approxChol wins) give a $56\%$ majority-class baseline, and the point-biserial
correlation of $\log\bar d$ with a LAMG+ win is $r=+0.68$ ($+0.72$ at $\varepsilon=10^{-4}$). The
mean degree thus predicts the faster solver a priori at either accuracy.}\label{fig:predictor}
\end{figure}

We construct a simple predictor of a LAMG+ win from a single statistic computable before either solver
runs: the mean degree $\bar d=2m/n$. Treating ``LAMG+ wins iff $\bar d\ge30$'' as a binary classifier
over the $201$ large graphs, its accuracy peaks at $91\%$ near $\bar d\approx30$
(Figure~\ref{fig:predictor}), and the point-biserial correlation of $\log\bar d$ with a LAMG+ win is
$r=+0.68$. The mechanism is two cost axes keyed to properties that anti-correlate across this corpus:
approximate-Cholesky's fill grows with elimination-clique size, hence with degree, while LAMG+'s
operator complexity grows with poor separability (low-diameter, small-world structure)---and here the
well-separated graphs are the high-degree ones, so $\bar d$ proxies both. High $\bar d$ favors LAMG+
and low $\bar d$ favors sampled elimination; the residual $9\%$ are exactly the graphs where the two
axes disagree---high-degree but small-world, or low-degree but cleanly separable.

\paragraph{Does a richer model help?} We tested a learned selector to try to explain the residual $9\%$: a gradient-boosted decision tree (a histogram model of the LightGBM class~\cite{lightgbm}) trained on $13$
further a-priori features of the degree distribution---median, maximum, $90$th/$99$th percentiles,
standard deviation, coefficient of variation, skewness, hub ratio $\max/\bar d$, and the leaf- and
high-degree fractions, all computable before any solve. Under repeated stratified $5$-fold
cross-validation ($20$ repeats; a single $80/20$ split is too small to be stable at $n{=}201$) it
attains $89\pm4\%$---no better than the one-line $\bar d$ threshold ($91\%$); permutation importance concentrates almost entirely on $\bar d$, every other feature contributing $\le0.01$. The degree distribution thus carries no
a-priori signal beyond its mean; the residual $9\%$ is governed by separability, a global property no
cheap local statistic exposes, so the principled one-statistic rule is also the near-optimal one. A
black-box library can apply it directly---dispatch to LAMG+ when $\bar d\gtrsim30$ and to
approximate-Cholesky otherwise---to pick the faster solver per input at no measurement cost.

\subsection{Robustness across the GKS test families}\label{sec:families}
The crossover above is a speed result; robustness---converging at all---is a separate axis, and the
one on which graph-Laplacian solvers most often differ. Gao, Kyng, and Spielman (GKS) \cite{gks2023}
report AC as the only solver in their study (alongside CMG, hypre/BoomerAMG, PETSc, and ICC) to
converge across their full family set---random \texttt{chimera} graphs and randomly-weighted
variants, SDDM ``chimeras'', adversarial ``Sachdeva stars'', an SPE benchmark, and
uniform/anisotropic/high-contrast Poisson grids. That robustness is precisely why we take AC as the
near-linear baseline worth a same-machine comparison.

LAMG+ meets the same bar. With the standard right-hand side ($b=L\,x_{\mathrm{true}}$, zero-mean), it
converges at both $\varepsilon=10^{-8}$ and $\varepsilon=10^{-4}$ on every one of the benchmark's families generated by
GKS's own \texttt{Laplacians.jl}: \texttt{chimera}s and their weighted, uniform, and
semi-weighted variants ($96$ graphs, $n=10^4$ to $10^6$), SDDM chimeras ($12$), $3$-D Poisson grids
with uniform, anisotropic, and high-contrast coefficients ($12$), and Sachdeva-style stars ($3$). It
also converges on the benchmark's one \emph{external} family, the SPE fluid dataset: we reconstruct
the SPE10 reservoir model \cite{christie-blunt} directly from its public permeability field,
discretized by a two-point flux approximation into high-contrast (permeability ratio ${\sim}10^7$),
anisotropic graph Laplacians at three sizes ($n=2.6\times10^5$ to $1.1\times10^6$), each solved in
$5$ cycles. Extending \emph{past} the scale at which AC was benchmarked, it further solves $41$ large
SuiteSparse graphs up to \texttt{com-Orkut} ($1.2\times10^8$ edges, $2.4\times10^8$ nonzeros).
Together with the $1{,}711$-graph scaling set (\S\ref{sec:scaling}), every instance converges, in a
median of $4$ cycles. The original MATLAB LAMG was reported non-convergent on all matrices of this set
\cite{gks2023}. This is not a property of the algorithm or of its original code: running the
\emph{unmodified} LAMG~2.2.1---under GNU Octave (compatibility shims only; the multigrid algorithm,
setup, and MEX kernels unchanged), with \cite{gks2023}'s exact solver call
(\texttt{setup(`laplacian',$\cdot$)}, \texttt{errorReductionTol}$=10^{-8}$, default cycle cap)---on the
same families converges on every one (Table~\ref{tab:lamgorig}), as does
LAMG+, in as few or fewer cycles and to a strict $10^{-8}$ residual. This indicates that the reported failure was an evaluation artifact rather than an algorithmic
property; we do not attempt to diagnose the specific cause within the original benchmark harness.

\begin{table}[t]\centering
\caption{The original MATLAB LAMG 2.2.1 vs.\ LAMG+ on the families of Gao-Kyng-Spielman (GKS)
approximate-Cholesky benchmark \cite{gks2023}, \emph{same graph and right-hand side}. Both LAMG and LAMG+ converge for all instances. LAMG+ reaches a strict
$10^{-8}$ residual in as few or fewer cycles. ``cyc''~=~multigrid cycles; relres~=~$\|L\phi-b\|/\|b\|$ at
termination. (LAMG 2.2.1's \texttt{errorReductionTol} test targets the error norm, so its residual
lands within a small factor of $10^{-8}$, e.g.,\, $3.2\!\times\!10^{-8}$ for uniform-weight chimera.)}
\label{tab:lamgorig}
\footnotesize\setlength{\tabcolsep}{4pt}%
\begin{tabular}{l r r r r r r}
\toprule
 & & & \multicolumn{2}{c}{LAMG 2.2.1} & \multicolumn{2}{c}{LAMG+}\\
\cmidrule(lr){4-5}\cmidrule(lr){6-7}
Family & $n$ & $m$ & cyc & relres & cyc & relres\\
\midrule
chimera (unweighted) & 50,000 & 173,059 & 10 & $1.8\!\times\!10^{-9}$ & 5 & $8.5\!\times\!10^{-9}$\\
chimera (uniform wt) & 50,000 & 99,552 & 9 & $3.2\!\times\!10^{-8}$ & 5 & $4.9\!\times\!10^{-10}$\\
chimera (weighted) & 50,000 & 76,068 & 8 & $1.0\!\times\!10^{-9}$ & 4 & $1.6\!\times\!10^{-9}$\\
chimera (SDDM) & 50,000 & 173,059 & 10 & $2.7\!\times\!10^{-9}$ & 5 & $6.6\!\times\!10^{-9}$\\
Sachdeva star & 10,001 & 495,100 & 4 & $2.7\!\times\!10^{-10}$ & 3 & $4.6\!\times\!10^{-10}$\\
grid 2D (uniform) & 65,536 & 130,560 & 11 & $3.8\!\times\!10^{-9}$ & 5 & $6.4\!\times\!10^{-10}$\\
grid 2D (anisotropic) & 65,536 & 130,560 & 19 & $2.3\!\times\!10^{-9}$ & 6 & $7.8\!\times\!10^{-9}$\\
grid 2D (high-contrast) & 65,536 & 130,560 & 10 & $6.0\!\times\!10^{-9}$ & 4 & $6.6\!\times\!10^{-9}$\\
grid 3D (uniform) & 64,000 & 187,200 & 10 & $8.5\!\times\!10^{-9}$ & 9 & $6.1\!\times\!10^{-9}$\\
grid 3D (anisotropic) & 64,000 & 187,200 & 16 & $8.1\!\times\!10^{-9}$ & 11 & $2.2\!\times\!10^{-9}$\\
grid 3D (high-contrast) & 64,000 & 187,200 & 12 & $7.1\!\times\!10^{-9}$ & 12 & $3.7\!\times\!10^{-9}$\\
SPE10 (reservoir) & 264,000 & 773,200 & 11 & $4.8\!\times\!10^{-9}$ & 8 & $6.9\!\times\!10^{-9}$\\
\bottomrule
\end{tabular}
\end{table}

Our first-hand results corroborate the central finding of \cite{gks2023}: BoomerAMG, CMG, pyAMG,
and PETSc GAMG each fail to converge or slow by an order of magnitude on chimera, SDDM-chimera, or
adversarial-star instances (Table~\ref{tab:classcomp}); the present comparison adds LAMG+ as the
only other solver, alongside AC, that is robust across all twelve families. On these well-structured families approxChol is faster---bounded degrees and separators
suit its sampled elimination---while LAMG+'s advantage is on large poorly-separable graphs.

\subsection{Lower accuracy shifts the balance toward LAMG+}\label{sec:tol4}
Many applications need only a few digits. At $\varepsilon=10^{-4}$ the LAMG+ solve collapses
to a \emph{single} cycle. This helps LAMG+ across the board
(Table~\ref{tab:compete}, bottom block) for a concrete reason: LAMG+ spends the larger share of its
time in the \emph{solve} (a median $68\%$ of its total on the low-diameter class, versus $37\%$ for
approximate-Cholesky)---a few cycles over its higher-complexity hierarchy outweigh the lean
PCG---so shedding it pulls every graph toward its setup-only ratio,
where LAMG+ stands relatively better, even on the graphs it still loses. Against the approxChol the aggregate geometric mean rises from
$0.99\times$ at $\varepsilon=10^{-8}$ (parity) to $1.20\times$ at $\varepsilon=10^{-4}$, and its win rate
over the $201$ large graphs rises from $56\%$ to $64\%$. Against the robust AC the lead widens from
$2.20\times$ to $2.91\times$ ($76\%\!\to\!84\%$ of head-to-heads). The gain is broad but does not
overturn the low-diameter classes---there it only narrows LAMG+'s deficit; the decisive flips are the
moderate-degree graphs near the crossover. The degree-threshold predictor tracks this shift: its optimal cutoff reduces from $\bar d\approx30$ at $\varepsilon=10^{-8}$ to $\approx13$ at $\varepsilon=10^{-4}$
(Figure~\ref{fig:predictor}).

\subsection{Anisotropy: the veto and the caliber-2 extension}\label{sec:aniso-results}
Grid-aligned anisotropy is baseline (caliber-1) LAMG+'s worst class. On a $64\times64$ grid with
strong weight $1$ and weak weight $\epsilon$ the full solver converges at a bounded but
non-textbook $\acf\approx0.27$--$0.39$ for $\epsilon\le10^{-1}$ (vs.\ $0.01$ isotropic); the
caliber-2 gate (\S\ref{sec:cal2}) recovers the ideal rate $\acf=0.007/0.05/0.10$ at
$\epsilon=10^{-1}/10^{-2}/10^{-4}$, at \emph{lower} operator complexity (sharper interpolation
buys fewer levels). On larger grids ($256^2$, $384^2$) a second, rarer failure appears:
the affinity inversion of \S\ref{sec:soc} produces a single bad merge that stalls even caliber-2
at $\acf\approx0.99$. The strength-of-connection veto removes it, and the two refinements together
hold $\acf\approx0.11$ across all sizes, seeds, and $\epsilon\le10^{-1}$ (worst $0.117$ over $24$
cases). On the SuiteSparse stiffness matrices the pair converts $100$-cycle near-failures into
$5$--$14$ cycle solves. This does not make LAMG+ the fastest solver on this class---these
matrices are small and well-separated, so approxChol and CHOLMOD stay ahead---so the
value here is robustness, closing LAMG's documented outlier tail. The
mechanism, the per-level two-level factors, and interpolation weight derivation are in App.~\ref{sec:aniso}.

\subsection{Ablation Study}\label{sec:ablation}
Each choice---the reproduced LAMG choices, the cycle-index calibration, and the two refinements
of \S\ref{sec:choices}---is load-bearing only on the class that stresses it. Table~\ref{tab:ablate}
lists their individual effects. The dominant setup cost is the test-vector smoothing, fixed at the cheap
$K{=}4$ everywhere (\S\ref{sec:calib}); the anisotropy fix comes from caliber-2, not from more test
vectors (Table~\ref{tab:adaptive}). The veto and
caliber-2 are \emph{complementary} on their target class: on the $256^2$, $\epsilon=10^{-4}$ grid
the veto alone reaches its caliber-1 ceiling ($\acf\approx0.46$), caliber-2 alone still stalls at
$0.99$ (the affinity inversion survives), and only both together give $0.11$.

\begin{table}[t]\centering\footnotesize\setlength{\tabcolsep}{4pt}
\caption{When each choice is necessary, and its effect when wrong (negligible elsewhere).
``Stresses'' is the graph class on which the choice is load-bearing. Top: LAMG choices LAMG+
reproduces (\S\S\ref{sec:alg},\,\ref{sec:calib}); bottom: the two anisotropy refinements.}\label{tab:ablate}
\begin{tabular}{>{\raggedright\arraybackslash}p{2.9cm}
               >{\raggedright\arraybackslash}p{2.9cm}
               >{\raggedright\arraybackslash}p{5.0cm}}
\toprule
Choice & Stresses & Effect if wrong\\
\midrule
no size cap & scale-free / web / social & web $\acf\ 0.002\!\to\!0.97$ (dense fill, stall)\\
singletons & anisotropic FE & \texttt{bodyy5} diverges ($\to0.99$)\\
directional guard & structured & over-aggregation; rarely fatal alone\\
\midrule
strength-of-connection veto & grid-aligned anisotropy & rare affinity-inversion stall ($\to0.99$)\\
caliber-2 (default on) & grid-aligned anisotropy & ceiling $0.46\!\to\!0.11$; $100$-cyc near-failures $\to5$--$14$\\
\bottomrule
\end{tabular}
\end{table}

\section{Conclusion}\label{sec:conclusion}
LAMG+ is a lean, parameter-free, empirically linear-time graph-Laplacian solver. It reproduces LAMG's choices, plus two targeted, \emph{local} refinements (\S\ref{sec:choices})---a
strength-of-connection veto and a selective caliber-2 interpolation. LAMG+ works for structured, web, social, citation, collaboration, road, and
grid-aligned anisotropic graphs with no per-graph tuning, and on every family of GKS's
approximate-Cholesky benchmark \cite{gks2023} without fail. LAMG+ and approximate-Cholesky are peers with complementary strengths---approximate-Cholesky faster on
social and citation graphs, LAMG+ the fastest and most memory-frugal solver on
finite-element/structural matrices. Beyond the linear solver, the multilevel hierarchy could be used
for related tasks such as eigenproblems \cite[\S6.4]{lamg} and nonlinear systems via the Full
Approximation Scheme \cite[\S8.1]{guide}.

\paragraph{On the feasibility of a parallel LAMG+} While the empirical results presented here are
single-threaded, the lean algebraic modifications introduced preserve the theoretical concurrency
bounds of the original solver. Both the setup and the cycle reduce to scalable sparse primitives---sparse
matrix--vector products, maximal independent set, segmented reductions---already mapped to GPUs
\cite{bell-dalton-olson} and scaled to $10^5$ cores \cite{baker-scaling,amgx}. The one inherently
sequential kernel is the Gauss--Seidel smoother; the standard substitutes are polynomial smoothers
\cite{adams-poly} or a block-Jacobi Gauss--Seidel that freezes the cross-block coupling and, in our
measurements, retains most of serial GS's smoothing power once the blocks are wide. The open obstacle
is distributed coarsening: a low-diameter Laplacian must be repartitioned per coarse level to bound
inter-processor communication \cite{metis}, and whether LAMG+'s purely local affinity and caliber-2
tests preserve its lean operator complexity under a communication-minimizing partition is unresolved.
A no-fill aggregation hierarchy is, if anything, friendlier to distribution than a sparsified Cholesky
factor whose elimination order is global.

\appendix
\section{Grid-aligned anisotropy: caliber-2 fix}\label{sec:aniso}

\begin{figure}[t]\centering
\includegraphics[width=0.82\textwidth]{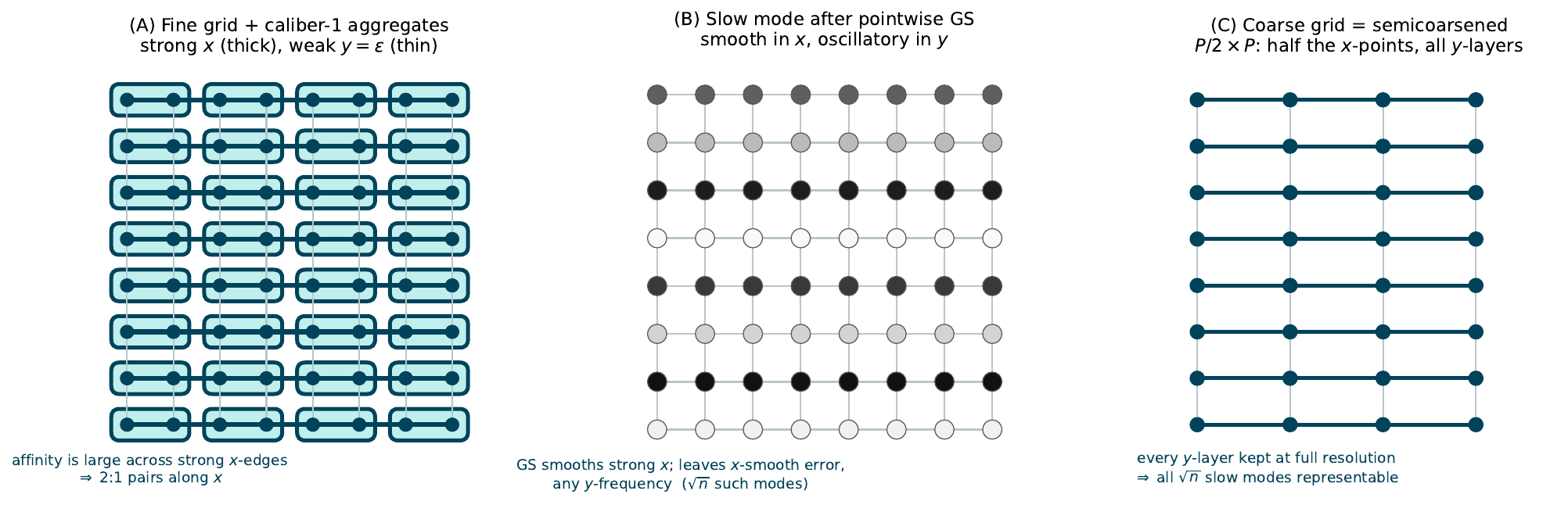}
\caption{Why caliber-2 is a one-dimensional fix. (A)~Affinity merges along the strong
direction ($2{:}1$ $x$-pairs)---a semicoarsening. (B)~The error pointwise GS leaves is smooth
in $x$, oscillatory in $y$. (C)~The semicoarsened grid keeps every $y$-layer; caliber-1
interpolates the smooth-$x$ mode across each $x$-pair by a constant ($q\approx2$,
$\rho\approx0.5$), caliber-2 by a line ($q\to1$, $\rho\approx0.12$).}\label{fig:semicoarsen}
\end{figure}

\begin{figure}[t]\centering
\includegraphics[width=0.82\textwidth]{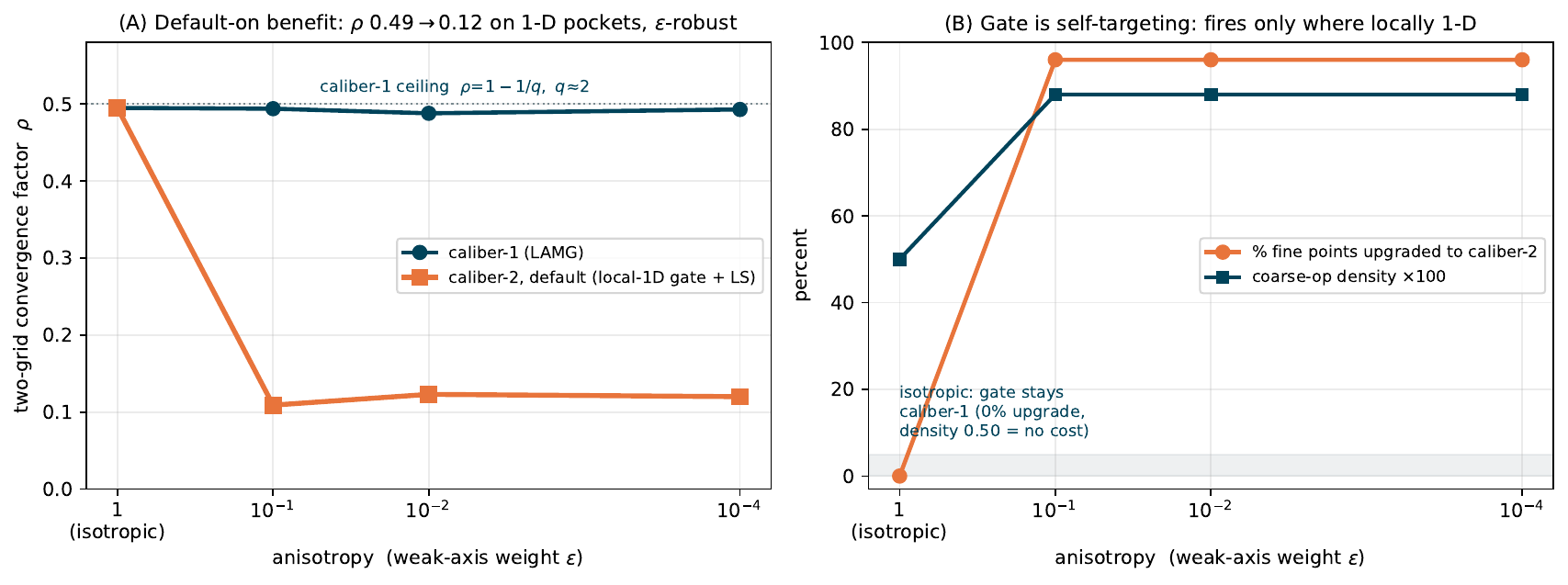}
\caption{The caliber-2 fix (\S\ref{sec:cal2}, on by default), two-grid factor on an anisotropic
grid. (A)~On locally one-dimensional pockets the second parent drops $\rho$ from the
caliber-1 ceiling $\approx0.49$ to $\approx0.12$, $\epsilon$-robust. (B)~The gate is
self-targeting: $\approx\!96\%$ of fine nodes upgrade where the neighborhood is $1$-D and
$\approx\!0\%$ on the isotropic grid.}\label{fig:caliber2}
\end{figure}

On a five-point grid with strong (horizontal) weight $1$ and weak (vertical) weight $\epsilon$,
the headline inverts the naive expectation: an isotropic grid is easy ($\acf=0.01$),
\emph{heterogeneous} weights (each edge an independent random factor over three decades) are
also easy ($\acf=0.006$), but \emph{grid-aligned} anisotropy is the hard case
($\acf\approx0.34$--$0.39$ at $\epsilon\le10^{-2}$). The mechanism and the cure are both in the
LAMG design. Affinity \emph{does} semicoarsen---it merges along the strong direction, the right
\emph{kind} of coarsening given a pointwise GS smoother. The caliber-1 piecewise-constant interpolation recovers half the magnitude of a
smooth-in-strong mode, inflating the coarse-operator energy by $q\approx2$ and pinning the
\emph{two-level} factor at $\rho\approx1-1/q\approx0.5$ \cite[\S3.4.3]{lamg}. In a multilevel setting, this slowness compounds across levels: a flat V-cycle ($\gamma=1$, no recombination) sits at $\acf\approx0.77$. Increasing $\gamma$ to $1.5$ and employing iterate recombination at every level \cite[\S3.4.5]{lamg}, we can bound the ACF at $0.39$ (Table~\ref{tab:aniso}). We can either accept this limited convergence factor (LAMG) or increase the interpolation caliber to improve it (LAMG+).

\subsection{A two-grid local Fourier analysis}\label{sec:lfa}
The ceiling $\rho\approx\tfrac12$ and its removal by
caliber-2 follow from a standard two-grid local Fourier analysis \cite[\S4.1]{guide}. Consider the
constant-coefficient infinite grid (strong weight $1$ along $x$, weak weight $\epsilon$ along $y$ on
$\mathbb Z^2$); the Laplacian acts on the error mode
$\varphi_\theta(\mathbf m)=e^{\mathrm{i}\,\theta\cdot\mathbf m}$,
$\theta=(\theta_x,\theta_y)\in(-\pi,\pi]^2$. Abbreviate $c=\cos^2\tfrac{\theta_x}{2}$ and
$s=\sin^2\tfrac{\theta_x}{2}$ ($c{+}s=1$); then
\begin{equation}\label{eq:Lsymbol}
\hat L(\theta)=4s+4\epsilon\sin^2\tfrac{\theta_y}{2}.
\end{equation}
Since $\epsilon\ll1$, $\hat L(\theta)$ is small exactly when $\theta_x\approx0$ for any $\theta_y$ (the
algebraically smooth error); affinity therefore semicoarsens $2{:}1$ in $x$ and keeps every $y$-line
(Fig.~\ref{fig:semicoarsen}).

\paragraph{Two-level cycle symbol}
Coarsening in $x$ at fixed $\theta_y$, the $2{:}1$ step folds each low mode $\theta$ ($\theta_x\in(-\tfrac\pi2,\tfrac\pi2]$) with its $x$-harmonic $\theta^{\ast}=(\theta_x{+}\pi,\,\theta_y)$ (same $\theta_y$); the pair spans a subspace invariant under the two-grid operator, on which the operator is the $2\times2$ matrix
\begin{equation}\label{eq:tgsymbol}
\hat T(\theta)=\hat S(\theta)^{\nu}\Big[\,I-\hat P\,(\hat P^{\ast}\hat L\hat P)^{-1}\hat P^{\ast}\hat L\,\Big],\qquad
\hat L=\mathrm{diag}\big(\hat L(\theta),\hat L(\theta^{\ast})\big)=:\mathrm{diag}(\ell_1,\ell_2),
\end{equation}
with the $2\times2$ diagonal smoother symbol $\hat S(\theta)=\mathrm{diag}\big(\tilde S(\theta),\,\tilde S(\theta^{\ast})\big)$,
$\nu=\nu_{\mathrm{pre}}{+}\nu_{\mathrm{post}}=3$ Gauss--Seidel sweeps, $\hat P\in\mathbb C^{2\times1}$ the
interpolation symbol, and the bracketed coarse-grid correction the $\hat L$-orthogonal projector off
$\mathrm{range}\,\hat P$.

\emph{Smoother.} $\tilde S$ is the scalar forward Gauss--Seidel symbol $\tilde S(\theta)=1-\hat L(\theta)/\hat M(\theta)$,
where $\hat M(\theta)=(2{+}2\epsilon)-e^{-\mathrm{i}\theta_x}-\epsilon\,e^{-\mathrm{i}\theta_y}$ is the
lower-triangular (already-updated-neighbor) part of $L$.

\emph{Interpolation.} Interpolating the fine point midway between two coarse parents has the symbol (up to a
unimodular phase that cancels in \eqref{eq:tgsymbol})
\begin{equation}\label{eq:Psymbols}
\hat P_{\mathrm{c1}}\propto\begin{pmatrix}\cos\tfrac{\theta_x}{2}\\[2pt]\mathrm{i}\,\sin\tfrac{\theta_x}{2}\end{pmatrix}
\ \ (\text{caliber-1, piecewise constant}),\qquad
\hat P_{\mathrm{c2}}=\begin{pmatrix}c\\ s\end{pmatrix}\ \ (\text{caliber-2, linear}).
\end{equation}
The decisive difference is the \emph{harmonic} (second) component as $\theta_x\to0$: it is $O(\theta_x)$ for caliber-1 but $O(\theta_x^2)$ for caliber-2. Caliber-2 injects asymptotically no energy into the unresolved harmonic of a smooth mode---the order-$2$ interpolation requirement of \cite[\S4.3]{guide}.

\emph{Energy factor.} The coarse-grid correction in \eqref{eq:tgsymbol} is a rank-one projector, so $\det\hat T=0$
and the spectral radius equals the trace:
\begin{equation}\label{eq:rhotrace}
\rho\big(\hat T(\theta)\big)=\big|\,\tilde S(\theta)^{\nu}(1-\beta)+\tilde S(\theta^{\ast})^{\nu}\,\beta\,\big|,
\qquad
\beta=\frac{|\hat P_1|^2\,\ell_1}{|\hat P_1|^2\,\ell_1+|\hat P_2|^2\,\ell_2}\in[0,1],
\end{equation}
where $\beta$ is the fraction of the interpolant's $\hat L$-energy carried by the smooth harmonic and
$q:=1/\beta$ is the local-mode form of the per-node energy ratio \eqref{eq:guard}. Substituting
\eqref{eq:Psymbols} and letting $\epsilon\to0$,
\begin{equation}\label{eq:betas}
\beta_{\mathrm{c1}}=\tfrac12\ \ (q=2)\quad\text{for every }\theta_x,
\qquad
\beta_{\mathrm{c2}}=c\to1\ \ (q\to1)\ \text{ as }\theta_x\to0 :
\end{equation}
caliber-1 commits exactly half of every interpolant's energy to the wrong harmonic, caliber-2 a
vanishing fraction on the smooth modes.

\emph{The caliber-1 ceiling.} The supremum of \eqref{eq:rhotrace} is attained as $\theta_x\to0$, on modes
\emph{smooth in $x$}, where any local smoother is powerless ($|\tilde S(\theta)|\to1$) while the
$x$-harmonic is fully smoothed ($\tilde S(\theta^{\ast})\to0$). There $\rho\to1-\beta$, so
$\beta_{\mathrm{c1}}=\tfrac12$ pins
\begin{equation}\label{eq:c1ceiling}
\rho_{\mathrm{c1}}\ \xrightarrow{\theta_x\to0}\ 1-\beta_{\mathrm{c1}}\;=\;1-\tfrac1q\;=\;\tfrac12\qquad(q=2).
\end{equation}
This ACF does not improve with more relaxation sweeps per cycle.

\emph{Caliber-2 removes it.} On those same smooth modes $\beta_{\mathrm{c2}}\to1$ ($q\to1$), so
$\rho\to1-\beta\to0$: the residual moves to the $x$-harmonic at $\theta_x=\tfrac\pi2$, a genuinely
high-frequency, smoothable mode. Evaluating \eqref{eq:rhotrace} under the same forward Gauss--Seidel
gives $\sup_\theta\rho_{\mathrm{c2}}\approx0.08$--$0.09$ for $\epsilon\le10^{-1}$---textbook two-grid efficiency.

Table~\ref{tab:aniso} compares caliber-1, the caliber-2 LFA prediction, and the measured caliber-2 two-grid factor against anisotropy strength. The LFA predicts the two-grid factor well; the small excess of the measured value is the least-squares weight's deviation from the ideal $\tfrac12$. Recombination adds a marginal improvement, and the full multilevel solver even beats the two-grid LFA rate for $\epsilon\ge10^{-2}$.

\begin{table}[h]\centering\footnotesize\setlength{\tabcolsep}{4pt}
\caption{Grid-aligned anisotropy on a $64\times64$ grid: asymptotic convergence factor by method.
Under caliber-1, \emph{V-cycle} $=$ flat $\gamma{=}1$ $(1,2)$-cycle and \emph{ML} $=$ the full solver
($\gamma{=}1.5$ $+$ recombination). Under caliber-2, \emph{LFA} is the local Fourier analysis
prediction, \emph{2-level} the measured two-grid factor without iterate recombination, \emph{$+$recomb}
the same with recombination, and \emph{ML} the full LAMG+ multilevel solver with the caliber-2
refinement. The two-grid columns (LFA, 2-level, $+$recomb) are measured on a periodic grid to match the
infinite-grid analysis; caliber-2 is inert at $\epsilon{=}1$ (isotropic, $0\%$ upgrade).}\label{tab:aniso}
\begin{tabular}{l cc cccc}
\toprule
 & \multicolumn{2}{c}{caliber-1} & \multicolumn{4}{c}{caliber-2}\\
\cmidrule(lr){2-3}\cmidrule(lr){4-7}
$\epsilon$ & V-cycle & ML & LFA & 2-level & $+$recomb & ML\\
\midrule
$1$       & 0.74 & 0.010 & ---  & ---  & ---  & \textbf{0.008}\\
$10^{-1}$ & 0.76 & 0.274 & 0.08 & 0.08 & 0.07 & \textbf{0.007}\\
$10^{-2}$ & 0.77 & 0.335 & 0.09 & 0.11 & 0.09 & \textbf{0.051}\\
$10^{-4}$ & 0.77 & 0.391 & 0.09 & 0.12 & 0.09 & \textbf{0.100}\\
\bottomrule
\end{tabular}
\end{table}

\subsection{Fitting the caliber-2 weight}
A fine node $u$ whose strong edges reach exactly two
coarse aggregates $a,b$ is interpolated by $\phi_u=w\,\phi_a+(1-w)\,\phi_b$; fixing the weights
to sum to $1$ keeps constants exact. The single weight $w$ is fit by least squares over the
\emph{same} $K$ test vectors used for affinity (no extra setup)---a two-parent specialization of the
least-squares interpolation of Bootstrap AMG \cite{brandt-bootstrap}: minimize
$\sum_k \big(x^{(k)}_u-[w\,x^{(k)}_a+(1-w)\,x^{(k)}_b]\big)^2$, whose closed form, with
$d_k=x^{(k)}_a-x^{(k)}_b$ and $e_k=x^{(k)}_u-x^{(k)}_b$, is
\begin{equation}\label{eq:wfit}
w^\star=\frac{\sum_k d_k\,e_k}{\sum_k d_k^2},
\end{equation}
the regression of $u$'s profile onto the strong chain's endpoints. We guard $w^\star\in[0,1]$
(a convex combination, so the Galerkin operator remains a Laplacian), snapping
near-boundary values to $0/1$; an out-of-range fit signals that $u$ is not truly one-dimensional
and is left caliber-1. The gate fires only when the two strong-reachable aggregates are distinct
and the fit is admissible---self-targeting at no parameter cost.

\subsection{The role of caliber-2}
The guard $Q$ of \eqref{eq:guard} already makes the solver robust in every
dimension ($\rho\lesssim1-1/Q$, never diverges). What it does not do is make $q$ \emph{small}
where anisotropy is persistent: grid-aligned anisotropy and $1$-D FE chains sit at the ceiling
$q\approx2$ at every level. Caliber-2 attacks that residue, and only where two conditions
coincide: the bound is tight (the strong coupling is a line, so caliber-1 across it is exactly
the half-recovered mode) and the cure is cheap (two collinear parents). Both hold iff a node is
locally one-dimensional. Hence anisotropic $2$-D and $3$-D \emph{rods} (a single strong
direction) become textbook, while $3$-D \emph{plates} (strong coupling in a plane) are locally
two-dimensional: the gate sees $>2$ strong aggregates and correctly declines, leaving this case to the guard.

\section{Linear-time evaluation of the energy guard}\label{sec:guardeval}
A direct evaluation of \eqref{eq:guard} is superlinear in node degree: recomputing, for each candidate
$(u,s)$ pair, the nodal energy by summing over $u$'s incident edges for each of $K$ test
vectors costs $\mathcal{O}(K\deg u)$ per pair, and a node with $D$ admissible neighbors is
tested $\mathcal{O}(D)$ times, so the work is $\mathcal{O}(K\deg^2)$---concentrated on the
high-degree hubs that dominate the total. The guard admits an $\mathcal{O}(K)$-per-candidate
form. In one $\mathcal{O}(Km)$ pass accumulate, per node $u$, the three moments
$s_0=\sum_v w_{uv}$, $s_1^{(k)}=\sum_v w_{uv}x^{(k)}_v$, $s_2^{(k)}=\sum_v w_{uv}(x^{(k)}_v)^2$.
Then both energies in \eqref{eq:guard} are closed forms: the relaxed (optimal) energy is
$\min_y E_u=\tfrac12\big(s_2^{(k)}-(s_1^{(k)})^2/s_0\big)$ and the forced energy is
$E_u(\cdot;x_s)=\tfrac12\big(x_s^2 s_0-2x_s s_1^{(k)}+s_2^{(k)}\big)$, so each candidate costs
$\mathcal{O}(K)$, independent of degree.

\end{document}